\documentclass{amsart}
\usepackage{amsmath,amssymb,graphicx}
\usepackage{mathrsfs}
\usepackage{pinlabel}

\newtheorem{theorem}{Theorem}[section]

\newtheorem{conjecture}[theorem]{Conjecture}

\newtheorem{definition}[theorem]{Definition}

\newcommand{\Z}{\mathbb{Z}}

\begin{document}
	
\title{COVID on trees and infinite grids}
\author[Barnett]{Andrea Barnett}
\author[Bond]{Robert Bond}
\author[Macias]{Anthony Macias}
\author[Mattman]{Thomas W.\ Mattman}
\email{TMattman@CSUChico.edu}
\author[Parnell]{Bill Parnell}
\author[Schoenfield]{Ely Schoenfield}
\address{Department of Mathematics and Statistics,
California State University, Chico,
Chico, CA 95929-0525}


\begin{abstract}
We use Hartnell's model for virus 
spread on a graph, also known as firefighting. For rooted trees, 
we propose an Unburning Algorithm, a type of greedy
algorithm starting from the leaves and working back towards the root.
We show that the algorithm saves at least half the vertices of the optimal
solution and that this is bound is sharp. We confirm a conjecture of Hartke
about integrality gaps when comparing linear and integer program solutions.
For general graphs, 
we propose a Containment Protocol, which looks ahead two time steps
to decide where to place vaccinations. We show that the protocol performs
near optimally on four well-studied infinite grids. The protocol is available
for any graph and we realize this flexibility by investigating an
infinite pentagonal graph.
 \end{abstract}
	

\maketitle

\section{Introduction}
We have recently experienced the COVID pandemic and know how a biological virus
can spread globally in a matter of months. Initially, vaccinations were in short supply, 
raising the question of how to quickly contain a virus with limited resources. 

In 1995, Hartnell~\cite{Ha} proposed a model for fire spread on a graph that has since been 
adapted to the spread of viruses or other contagions. 
We imagine a perfectly contagious virus spreading in discrete time steps through a population represented as
a graph. Each vertex is in one of three states: susceptible, infected, or vaccinated. Once a susceptible vertex is infected, the vertex remains 
infected and similarly for the vaccinated state. In other words, there is no cure and the only available response 
is that a finite number $n$ of susceptible vertices may be vaccinated at each time step. 

The outbreak begins with a single vertex infected. At each time step, we place $n$ vaccinations and then
the contagion spreads from the infected vertices to any neighbors that remain susceptible. The goal is to use the 
vaccinations to contain the outbreak with as few vertices as possible infected. The outbreak is contained when 
no infected vertex has a susceptible neighbor. A secondary and related goal
is to contain the outbreak in the least possible number of time steps.

In this paper we use this model to investigate virus spread on two types of graphs, rooted trees and infinite grids.
For rooted trees, the outbreak begins with the root vertex infected and we place one (ie $n = 1$) vaccine per time step. It is known
that the Greedy Algorithm will protect at least half as many vertices as an optimal solution. We propose the Unburning Algorithm,
which is a type of greedy algorithm, but working backwards from the leaves towards the root. We show that this algorithm
also performs quite well.
\begin{theorem}
\label{thmUBurn}
The Unburning Algorithm saves at least half as many vertices as saved by an optimal 
placement of vaccinations.
\end{theorem}

We prove this in Section~\ref{sec:UA} where we also present an infinite family of trees to show this bound is sharp.

Hartke~\cite{H2} made some modifications to an integer program of MacGillivray and Wang~\cite{MW} for virus spread
in a rooted tree. This led him to a linear program that is a good approximation to the optimal solution found by the
integer program. Indeed, Hartke conjectured that the linear program realizes the optimal solution for all trees of 12 or fewer vertices.
In Section~\ref{sec:H2} we verify this conjecture.

For infinite grids, we propose a Containment Protocol, which is a strategy for placing vaccinations based on looking ahead
two time steps. The Containment Protocol is quite robust and in Section~\ref{sec:IG} we show that it either matches or 
nearly matches the optimal solution for four well studied infinite grids, the hexagonal, square, triangular, and strong grids.
The Containment Protocol could be applied to any graph. As an example, in Section~\ref{sec:PG} we investigate 
an infinite pentagonal grid. We report on the best solutions found by the protocol that, we conjecture, are optimal.

\section{Trees}

In this section we investigate virus spread on a rooted tree. At time step  0, the outbreak begins with the root vertex infected. At each time
step we vaccinate one vertex. We will first contrast the Greedy Algorithm with our Unburning Algorithm. Then, in Section~\ref{sec:H2},
we verify a conjecture of Hartke~\cite{H2}.

\subsection{
\label{sec:UA}%
Greedy versus Unburning Algorithm}

In a rooted tree with root $v_0$, the {\em level $l(v)$} of  vertex $v$ is the distance or number of edges between $v_0$ and $v$: $l(v) = d(v_0,v)$.
Let $T_v$ denote the subtree rooted at vertex $v$.
Let $a_i$ be the vertex that is vaccinated at time step $i$. An optimal vaccination sequence $a_1, a_2, \ldots$ is one that saves
the maximum number of vertices. Hartnell and Li, and, independently, MacGillivray and Wang provided an important basic observation.

\begin{theorem}[\cite{HL,MW}] 
If $a_1, a_2, \ldots$ is optimal, then $l(a_i) = i$ for all $i$.
\end{theorem}

In the Greedy Algorithm, at each time step $i$, we vaccinate the 
vertex at level $i$ that is the root of the biggest remaining subtree. This protects the vertex
along with all the vertices in its subtree.
That is, if we define the weight of vertex $v$ as 
the number of vertices in the subtree rooted at $v$:
$\mbox{wt}(v) = |V(T_v)|$, then at time $i$, the Greedy Algorithm vaccinates the vertex $v$ that maximizes
$wt(v)$ among the remaining, unprotected vertices with $l(v) = i$. 

It is known that the Greedy Algorithm performs quite well, in general.

\begin{theorem}[Hartnell and Li~\cite{HL}] 
\label{thm:HL}
The Greedy Algorithm saves at least half as many vertices as saved by an optimal 
placement of vaccinations.
\end{theorem}

We next propose a new algorithm that we call the Unburning Algorithm. The name reflects the connection 
with firefighting. We imagine a fire sweeping through the tree and then attempt to `unburn' it starting from the leaves.

In the Unburning Algorithm, we work backwards from the maximum level $n$, choosing the sequence $u_n, u_{n-1}, u_{n-2}, \ldots$ in turn
with $l(u_i) = i$. To avoid some obvious poor choices, the weight of vertices at level $i$ are updated to reflect vaccinations that will
be placed below level $i$. Suppose $v$ has $l(v) = i$ and that $\{u_1^v, u_2^v, \ldots, u^{v}_{k_v}\}$ is the (possibly empty) set 
of vertices in $T_v$ included in the sequence $u_n, u_{n-1}, \ldots, u_{i+1}$. Let $U_v$ be the vertices in $T_v$ protected by 
those vaccinations:
$$ U_v = \bigcup_{k=1}^{k_v} V(T_{u_k^v}).$$
We define
$$\tilde{\mbox{wt}}(v) = |V(T_v)| - |U_v|.$$
and choose $u_i$ among the vertices $v$ with $l(v) = i$ to maximize $\tilde{\mbox{wt}}(v)$.

Our main observation is that, similar to the Greedy Algorithm, the Unburning Algorithm performs quite well, in general.

\setcounter{section}{1}
\setcounter{theorem}{0}

\begin{theorem}
The Unburning Algorithm saves at least half as many vertices as saved by an optimal 
placement of vaccinations.
\end{theorem}

\setcounter{section}{2}
\setcounter{theorem}{2}

\begin{proof}
We adapt Hartke's~\cite{H} charging argument for the Greedy Algorithm. Fix an optimal sequence of vaccinations $b_1, b_2, \ldots, b_k$ 
that saves the largest number of vertices and let $u_1, u_2, \ldots, u_{\ell}$ be the vertices selected by the Unburning Algorithm.
Whenever a vertex $u_i$ is in the subtree $T_{b_j}$ for some $b_j$ we will assign or `charge' $u_i$ to $b_j$. This means $i \geq j$.
We can assume that no $b_i$ is in the subtree $T_{b_j}$ with $j < i$ as 
we would then save the same number of vertices by removing $b_i$ from the sequence. Then each $u_i$ is charged to at most one
$b_j$.  Let $C_{b_j} = \{u_i \mid u_i \mbox{ is charged to } b_j\}$ denote the set of vertices charged to  $b_j$.

Let $S_o = \sum_{i=1}^k \mbox{wt}(b_i)$ denote the number of vertices saved by the optimal sequence.
Similarly, let $S_u = \sum_{i =1}^{\ell} \tilde{\mbox{wt}}(u_i)$ denote the number of vertices saved by the Unburning Algorithm. 
If, for some $i$, $\tilde{\mbox{wt}}(u_i) < \mbox{wt}(b_i)$, it must be that $C_{b_i} \neq \emptyset$ as otherwise the Unburning
Algorithm would have chosen $u_i = b_i$. Since the Unburning Algorithm chooses $u_i$ rather than $b_i$, we
know that the $\tilde{\mbox{wt}}(u_i)$ together with the weights of the $u_j$ assigned to $b_i$ is at least as big as
$\mbox{wt}(b_i)$. That is, 
$$\mbox{wt}(b_i) \leq \tilde{\mbox{wt}}(u_i) +  \sum_{u \in C_{b_i}} \tilde{\mbox{wt}}(u).$$
On the other hand, this inequality remains valid in case $\tilde{\mbox{wt}}(u_i) \geq \mbox{wt}(b_i)$.

Summing over the vertices in the optimal sequence, 
\begin{align*}
S_o &= \sum_{i=1}^k \mbox{wt}(b_i) \\
& \leq \sum_{i=1}^k \left(\tilde{\mbox{wt}}(u_i) +  \sum_{u \in C_{b_i}} \tilde{\mbox{wt}}(u) \right) .
\end{align*}

Fix an index $i_0$ in $1, \ldots, \ell$. There are at most two occurrences of $u_{i_0}$ in the sum. It will appear in $\tilde{\mbox{wt}}(u_i)$ when $i = i_0$.
Also, $u_{i_0}$ is the descendant of at most one $b_j$. So, $u_{i_0}$ may occur again in $C_{b_j}$ when $ i = j$. Since each $u_{i_0}$ occurs at 
most twice in the sum, we conclude that $S_o \leq 2S_u$, which means $S_u \geq \frac12 S_o$.
\end{proof}

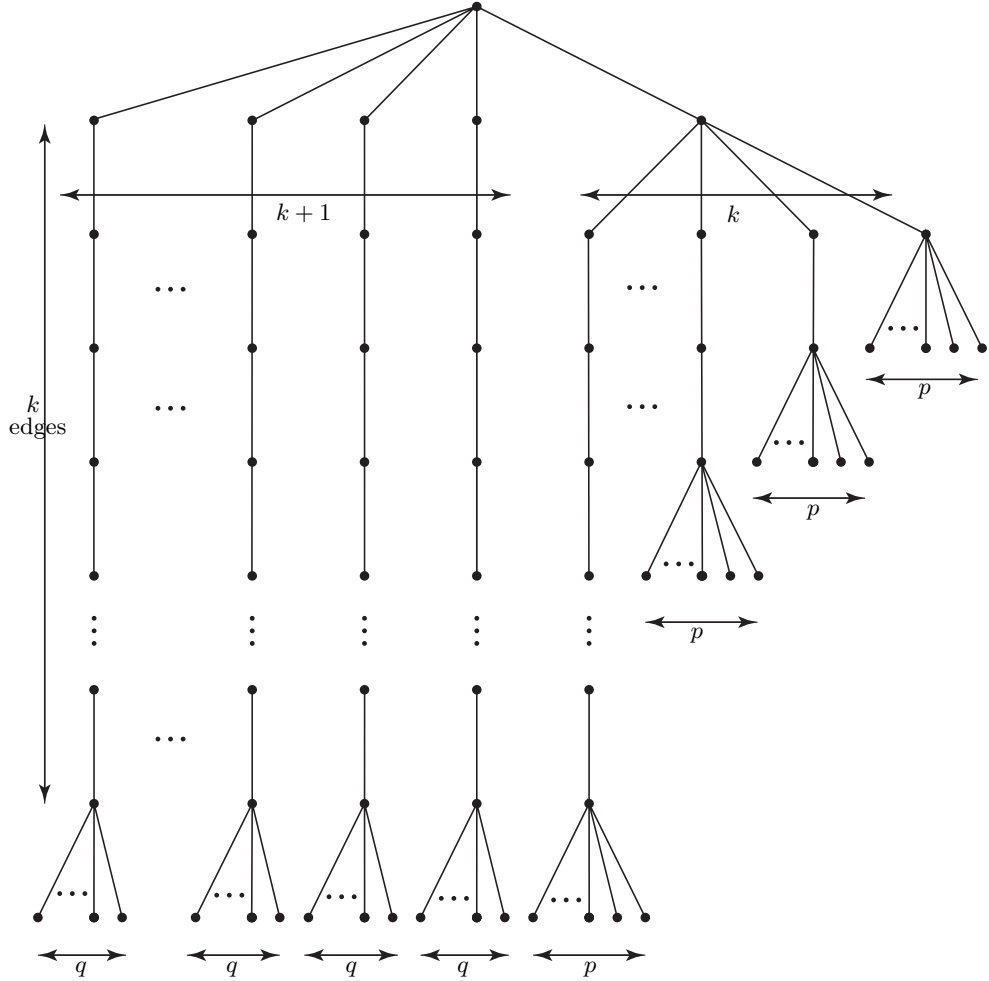
\begin{figure}[htb]
    \labellist                             
    \small\hair 2pt
    \pinlabel $q$ at 30 -5
    \pinlabel $q$ at 125 -5
    \pinlabel $q$ at 200 -5
    \pinlabel $q$ at 270 -5
    \pinlabel $p$ at 350 -5
    \pinlabel $p$ at 417 205
    \pinlabel $p$ at 490 282
    \pinlabel $p$ at 560 358
    \pinlabel $k$ at 0 350
    \pinlabel edges at 3 335
    \pinlabel $k+1$ at 170 470
    \pinlabel $k$ at 440 470
    \endlabellist         
\centering
\includegraphics[scale = 0.6]{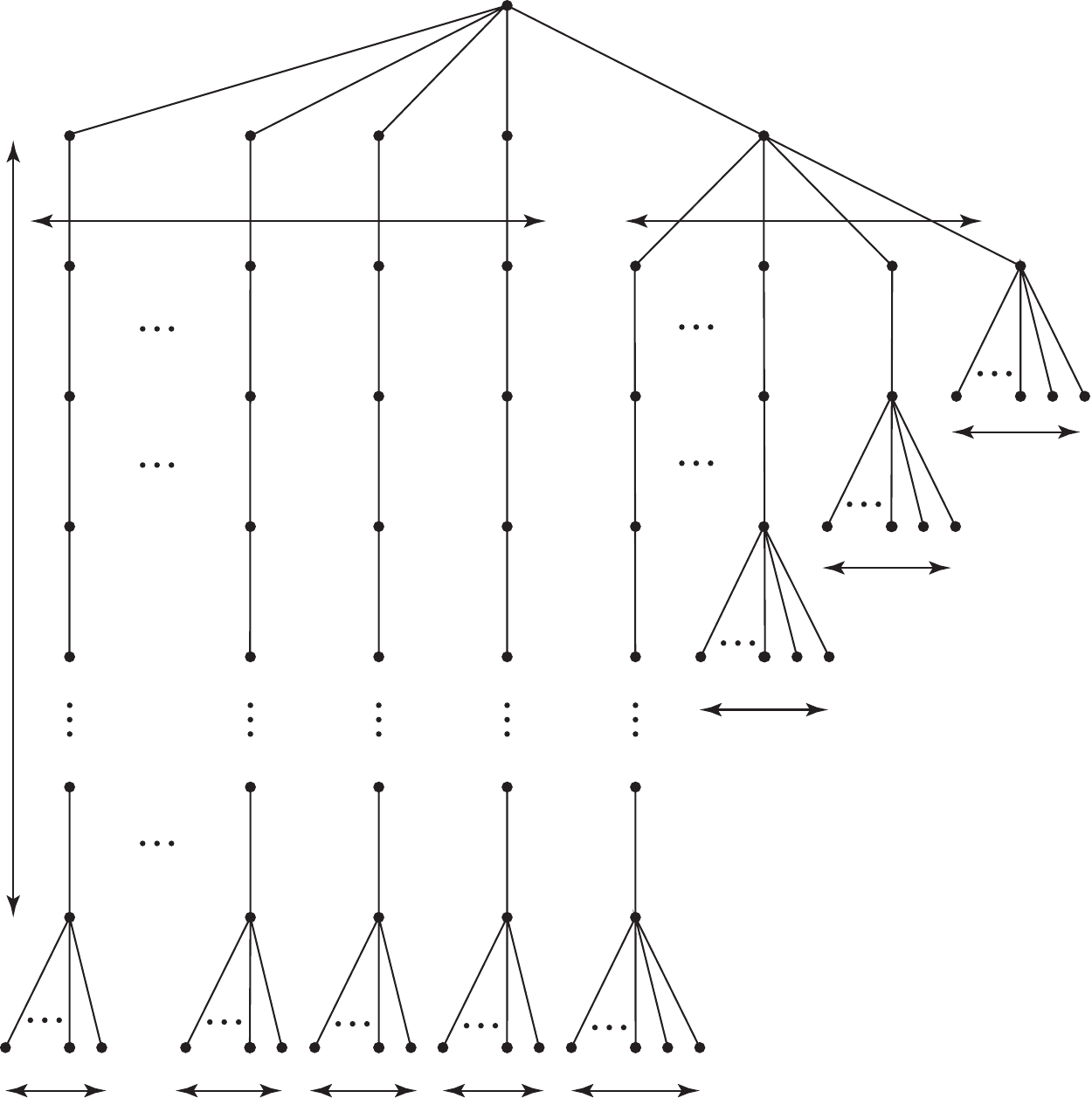}
\caption{The $(k,p,q)$ tree.}
\label{fig:kpq}
\end{figure}

We next show that this bound is tight by providing a family of trees for which the number of vertices saved by the Unburning Algorithm can be made 
arbitrarily close to half of those saved optimally. Let $(k,p,q)$ denote a tree of the type illustrated in Figure~\ref{fig:kpq}. The root vertex
has degree $k+2$. The $k+1$ subtrees to the left are identical, each consisting of a path of length $k$ with $q$ pendant vertices
as leaves. The remaining subtree to the right, includes $k$ subtrees, each with $p$ leaves at the ends of paths of length increasing
from $0$ to $k-1$. 

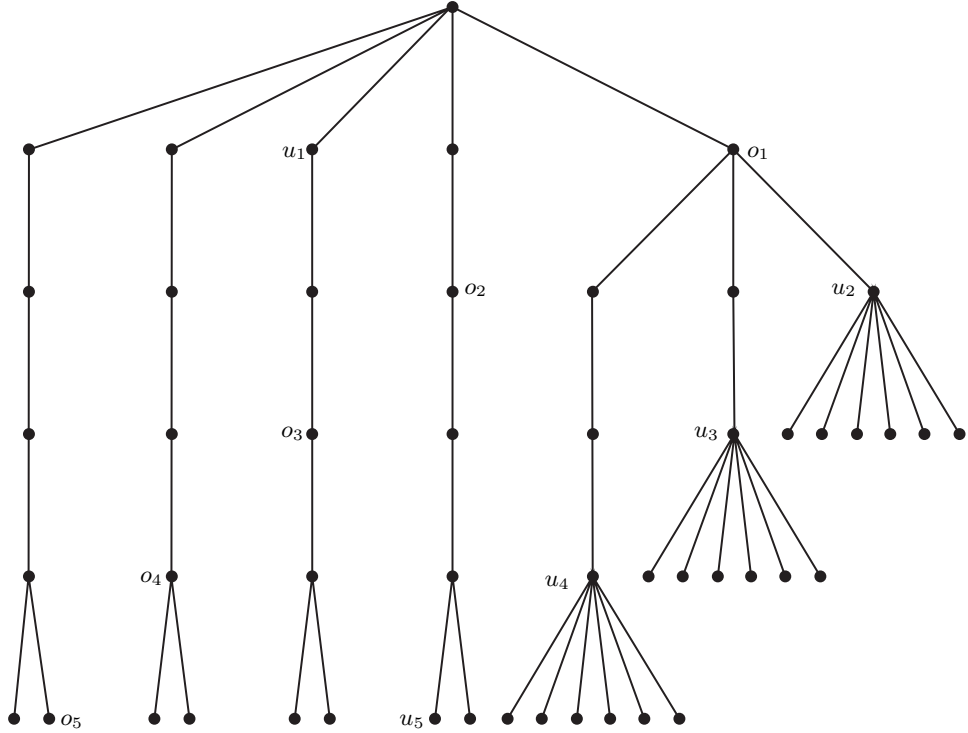
\begin{figure}[htb]
    \labellist                             
    \small\hair 2pt
    \pinlabel $o_1$ at 377 291
    \pinlabel $o_2$ at 235 222
    \pinlabel $o_3$ at 143 150
    \pinlabel $o_4$ at 72 77
    \pinlabel $o_5$ at 32 5
    \pinlabel $u_1$ at 144 291
    \pinlabel $u_2$ at 420 222
    \pinlabel $u_3$ at 351 150
    \pinlabel $u_4$ at 276 75
    \pinlabel $u_5$ at 203 5
    \endlabellist         
\centering
\includegraphics[scale = 0.75]{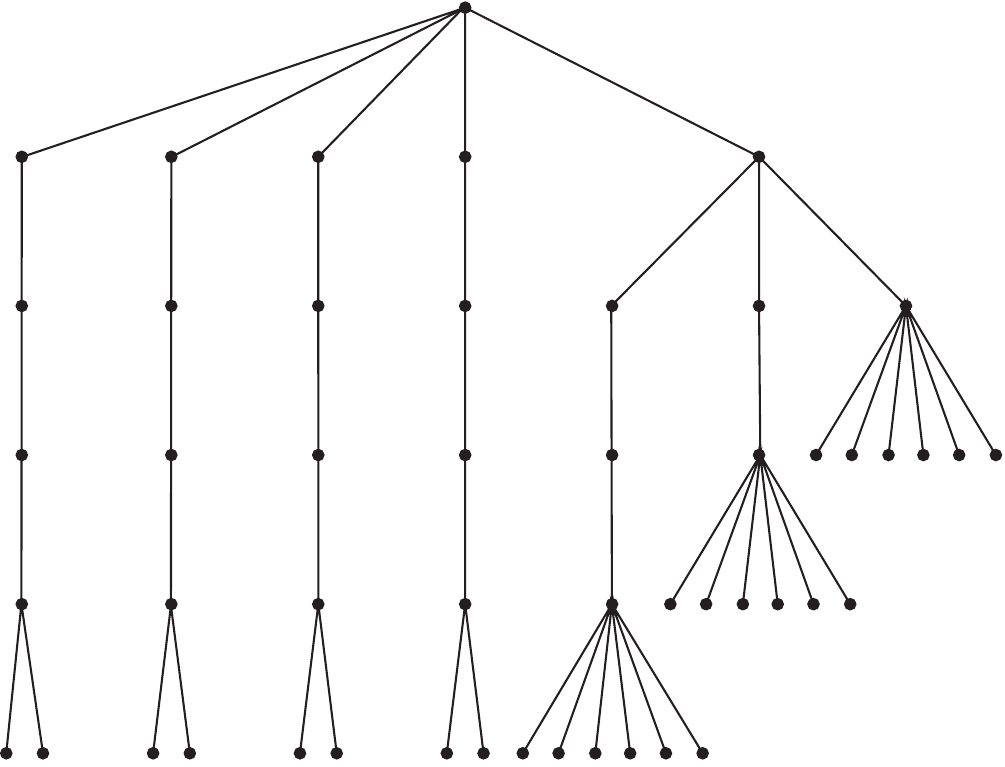}
\caption{The $(3,6,2)$ tree.}
\label{fig:362}
\end{figure}

For example, Figure~\ref{fig:362} illustrates the situation when $(k,p,q) = (3,6,2)$. For this graph, the Greedy Algorithm is also optimal. 
The vertices for this optimal solution $o_1, o_2, \ldots, o_5$ are as shown in the graph with total weight $S_o = 25 + 5+4+3+1 = 38$.
We have also labelled the vertices $u_1,u_2, \ldots, u_5$ of the vertices for the Unburning Algorithm with a total weight
$S_u = 6 + 7 + 7 + 7 + 1 = 28$. 

To generalize the $(3,6,2)$ example we will impose conditions on $k$, $p$, and $q$. 
Compare the vertices at level one. Note that 
$\tilde{\mbox{wt}}(o_1) = 4 < 6 = \mbox{wt}(u_1)$. In general, the subtree at the right in a $(k,p,q)$ tree will have
$\tilde{\mbox{wt}}(o_1) = k(k-1)/2 +1$ while the equivalent $k+1$ subtrees to the right have $\mbox{wt}(u_1) = k+q+1$. 
To maintain the relationship between these weights, we will require $k(k-1)/2 +1 < k+q+1$, or equivalently, $q > k(k-3)/2$.
Next look at the vertices at level two where $\mbox{wt}(u_2) = 7 > 5 = \mbox{wt}(o_2)$. In a $(k,p,q)$ tree we will have
$\mbox{wt}(u_2) = p+1$ and $\mbox{wt}(o_2) = k+q$. To maintain that relationship, we require that $p+1 > k + q$
or $p > k+q-1$. 

Assuming these conditions on $k$, $p$, and $q$, we can determine formulas for the number of vertices saved by the optimal
placement of vaccinations (which agrees with the Greedy Algorithm) and the Unburning Algorithm. For the optimal algorithm,
\begin{align*}
\mbox{wt}(o_1) &= 1 + (p+1) + (p+2) + \cdots + (p+k) \\
 & = 1 + kp + k(k+1)/2.
\end{align*}
The weights of $o_2, o_3, o_4, \ldots o_{k+1}$ go down by one at each step:
\begin{align*}
\sum_{i=2}^{k+1} \mbox{wt}(o_i) &= (q+k) + \cdots + (q+2) + (q+1) \\
&= kq + k(k+1)/2.
\end{align*}
Adding in $\mbox{wt}(o_{k+2}) = 1$, we have
$$S_o = 1 + kp + k(k+1)/2 + kq + k(k+1)/2 + 1 = k(k+p+q+1) + 2.$$

For the Unburning Algorithm, we have $\tilde{\mbox{wt}}(u_1) = \mbox{wt}(u_1) = k+1+q$, 
$\tilde{\mbox{wt}}(u_i) = \mbox{wt}(u_i) = p+1$ for $2 \leq i \leq k+1$, and
$\tilde{\mbox{wt}}(u_{k+1}) = \mbox{wt}(u_{k+1}) = 1$ so that
$$S_u = k+1+q + k(p+1) + 1 = k(p+2) + q + 2.$$

One way to satisfy our conditions on $k$, $p$, and $q$ is by setting $q = k(k-3)/2 + 1$ and 
$$p = k+q = k + k(k-3)/2+1 = k(k-1)/2 + 1.$$ With these choices, 
\begin{align*}
S_u/S_o &= \frac{k(k(k-1)/2 + 3) + k(k-3)/2 + 3}{k \left(k+k(k-1)/2+1 + k(k-3)/2+2 \right) +2} \\
 &= \frac{\frac12 k^3 + \frac32 k + 3}{k\left( k+3 + k(k-2) \right) + 2} \\
 &= \frac{\frac12 k^3 + \frac32 k + 3}{ k^3 -k^2 + 3k+2}.
 \end{align*}
Taking the limit as $k$ goes to infinity, the ratio becomes arbitrarily close to $\frac12$.
 
We remark that, for the $(k,p,q)$ trees, the Greedy Algorithm is optimal and the Unburning 
Algorithm does worse, tending to saving half as many vertices as $k$ increases.
In contrast, in his thesis, Hartke~\cite{H} shows that the bound of Theorem~\ref{thm:HL} is tight 
using a family of trees that are the right subtrees of the $(k,p,q)$ trees of Figure~\ref{fig:kpq}, that is
the analogs of the subtrees $T_{o_1}$ of Figure~\ref{fig:362}. For those trees, it is the Unburning Algorithm
that agrees with the optimal placement of vaccinations. Indeed, this was part of our initial 
motivation in developing the Unburning Algorithm. 

In general, the Greedy Algorithm and Unburning 
Algorithm are complementary in that if one performs poorly on a given tree, the other tends
to correct for that weakness. While the Greedy Algorithm and Unburning Algorithm share strengths
in that they are both straightforward to apply and perform at least half as well as the optimal
sequence of vaccinations, an even better general approach would be to apply both algorithms
and pick the one that saves more vertices.

\subsection{
\label{sec:H2}%
Integrality Gap Conjecture}

In this subsection we verify Conjecture 2.2 of Hartke~\cite{H2}, to which we refer the reader for further details.
MacGillvray and Wang~\cite{MW} formulated an integer program
to find an optimal sequence $a_1, a_2, \ldots$ for virus spread in a tree $T$. For each vertex $v$, define
$x(v) = 1$ if $v$ is vaccinated and $0$ otherwise. Recall that $l(v) = d(v_0,v)$ is the number of edges
between $v$ and the root vertex $v_0$ and $\mbox{wt}(v) = |V(T_v)|$ is
the number of vertices in the subtree rooted at $v$. If vertex $u$ is in $T_v$, we will write $v \succeq u$.

With this notation, we can write MacGillvray and Wang's
integer program as follows.
$$ \mbox{maximize } m = \sum_{v \in V(T)} wt(v)x(v)$$
\begin{align}
\mbox{subject to: }  \sum_{l(v) =  L} x(v) \leq 1, & \mbox{ for each level }L,\\
 \sum_{v \succeq u} x(v) \leq 1, & \mbox{ for each leaf vertex } u, \mbox{ and } \\
 x(v) \in \{0,1\}, & \mbox{ for each vertex } v.
\end{align}

By relaxing constraint (3) that $x(v)$ be boolean, we instead obtain a linear program (LP) whose maximum $m^*$ is an upper
bound for the maximum $m$ of the integer program (IP).  We say there is an ``integrality gap'' when $m^* > m$. Since
a linear program is much more efficient computationally, it is useful to understand the occurrence and size of integrality gaps.

Hartke~\cite{H2} attempts to reduce the occurrence of integrality gaps by adding an additional linear constraint, constraint (6).
\begin{align*}
(6) \hspace{0.5 in}  \sum_{v \succeq u} x(v) + \sum_{u \succeq v, l(v) = i} x(v) \leq 1, \mbox{ for each vertex } u \mbox{ and } i > l(u). 
\end{align*}

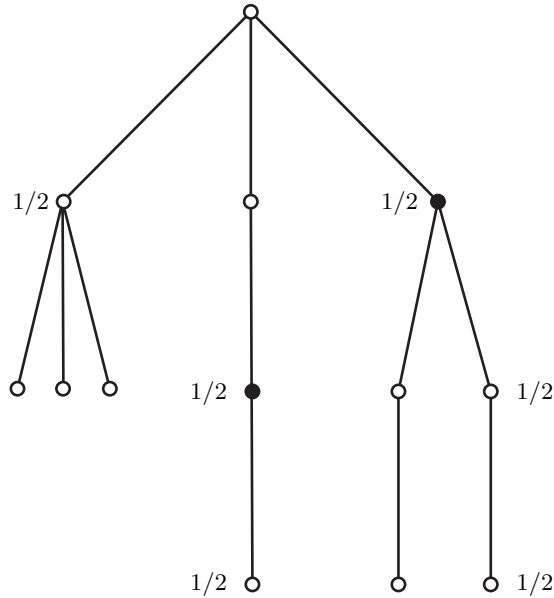
\begin{figure}[htb]
    \labellist                             
    \small\hair 2pt
    \pinlabel $1/2$ at 75 3
    \pinlabel $1/2$ at 198 3
    \pinlabel $1/2$ at 75 75
    \pinlabel $1/2$ at 198 75
    \pinlabel $1/2$ at 8 147
    \pinlabel $1/2$ at 147 147
    \endlabellist         
\centering
\includegraphics[scale = 1]{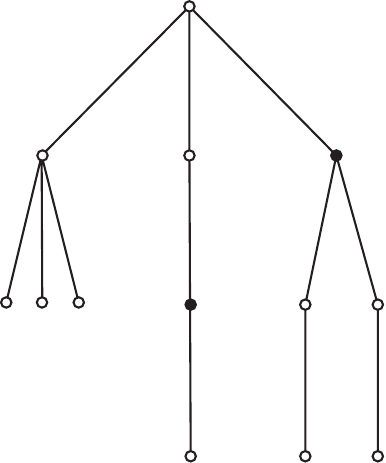}
\caption{A tree, due to Hartke~\cite{H2}, on 13 vertices with integrality gap when using constraint (6).}
\label{fig:H3}
\end{figure}

Although this improves the situation, there remain trees with an integrality gap including the one on 13 vertices shown in Figure~\ref{fig:H3}.
The two black vertices give the optimal (IP) solution with $m = 7$ vertices saved. The nonzero $x(v)$ appear in the figure and give 
the LP optimum $m^* = 7.5$. 

We can now state Hartke's conjecture.

\begin{conjecture}[Conjecture 2.2 of \cite{H2}] The tree in Figure~\ref{fig:H3} is the smallest tree such that the LP optimal when using constraint
(6) is not the IP optimal.
\end{conjecture}

\begin{figure}[htb]
    \labellist                             
    \small\hair 2pt
    \pinlabel $1/2$ at 143 3
    \pinlabel $1/2$ at 203 3
    \pinlabel $1/2$ at 68 76
    \pinlabel $1/2$ at 203 76
    \pinlabel $1/2$ at 17 147
    \pinlabel $1/2$ at 157 147
    \endlabellist         
\centering
\includegraphics[scale = 1]{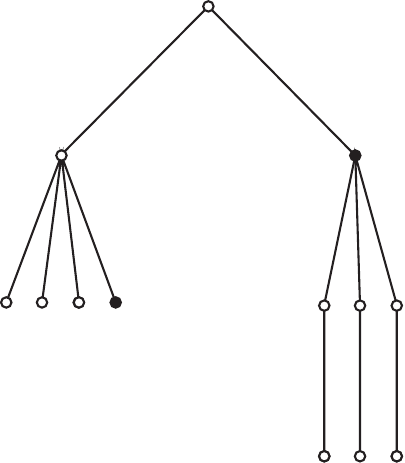}
\caption{A second tree on 13 vertices with integrality gap when using constraint (6).}
\label{fig:T13}
\end{figure}
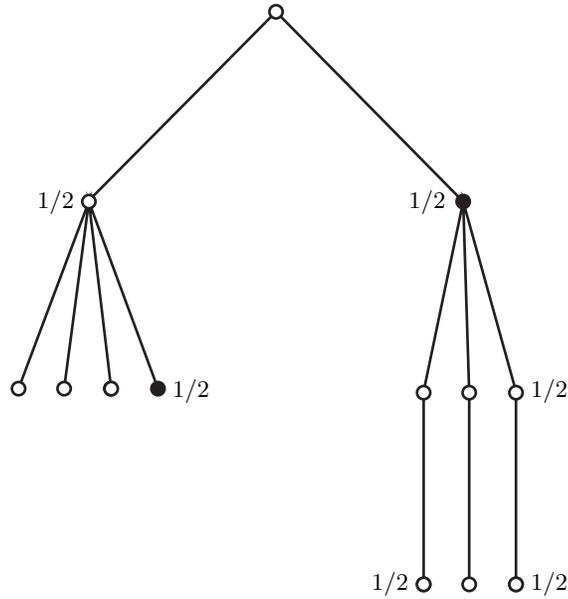

Here, ``smallest'' refers to the number of vertices. Hartke explained that he had verified the conjecture for trees of 11 or fewer vertices.
To complete the proof of the conjecture then, we used Sage~\cite{sage} to run through the 4766 trees on 12 vertices and check that none of those
have an integrality gap using constraint (6). We also checked trees on 13 and 14 vertices. Among the 12486 trees on 13 vertices, there are
exactly two with integrality gap. Figure~\ref{fig:T13} shows the second graph. The black vertices realize $m =   8$, while the given
$x(v)$ show that $m^* = 8.5$. We found that of the 32973 trees on 14 vertices, there are ten trees that have an integrality gap. 
Moreover, each of these ten is formed from one of the two examples on 13 vertices by adding a leaf and its edge. There are six
ways to add a leaf to the graph of Figure~\ref{fig:H3} and four ways for the graph of Figure~\ref{fig:T13}.
All 12 examples have $m^* - m = 0.5$.

\section{
\label{sec:IG}
Infinite regular graphs}

In this section we look at virus spread in four infinite regular graphs. Three correspond to tessellations of the plane by regular polygons:
the hexagonal, square, and triangular grids. A fourth model is the strong grid, which is regular of degree eight. 
We describe a Containment Protocol to bring an outbreak on a graph under control in a finite number of steps. We show that this strategy matches, or nearly
matches the best known results on these four grids.

We begin by describing the Containment Protocol. At each time step, we imagine allowing the virus to spread for two time steps. Based on this idea, we make
the following definitions.

\begin{definition}
At time step $T$, we denote by AIV, {\em already infected vertices}, the set of vertices that were infected at a time $t < T$, 
IV1, {\em infected vertices of type 1}, the vertices infected in time step $T$, and IV2, {\em infected vertices of type 2}, those infected in time step $T+1$ (assuming no vaccines are placed at time $T$). We refer to vertices where vaccinations were placed at times $t < T$ as PPV, {\em previously placed vaccines.}
For a $v \in $ PPV that is adjacent to an IV1 vertex, we say that $v$ is a {\em good PPV} if at most half of its neighbors are in AIV $\cup$ IV1. Otherwise, $v$ is a {\em bad PPV}.
\end{definition}

In the Containment Protocol, we place $n$ vaccines at each time step according to the following three criteria: 
0) each vaccine is placed at a vertex $v$ that is in neither AIV nor PPV,
1) each vaccine is placed adjacent to an AIV, and
2) each vaccine is at most distance two from either another vaccine placed at that time step, or a vaccine placed at an earlier time step, a PPV.

In rule CP5 below, we calculate the distance of the placement from the PPV vertices. For each $v$ in a candidate placement, let $d(v, \mbox{PPV})$ be
the distance to the PPV set. That is,
$$d(v, \mbox{PPV}) = \min_{u \in PPV} d(v,u),$$
where $d(v,u)$ is the number of edges in the shortest path from $v$ to $u$.

We choose among the allowed vaccine placements at a given time step according to the following five rules.
\begin{description}
\item[CP1] Choose a placement that minimizes $|\mbox{IV2}|$. If there is more than one placement that minimizes $| \mbox{IV2}|$ eliminate any placements
that do not and continue to the next rule.
\item[CP2] Choose arrangements that minimize the number of PPV adjacent to an IV2. If there is more than one that minimize this, continue with those
arrangements that minimize.
\item[CP3] If there is a placement that includes a good PPV, eliminate any placements that have no good PPV. If more than one
remains, continue to the next rule.
\item[CP4] If there remains a placement with no bad PPV, then eliminate any that have one.
If more than one placement remains, continue to the next rule.
\item[CP5] Use a placement that minimizes the total distance $\sum_{v} d(v, \mbox{PPV})$. If there is more than one placement that minimizes 
the total distance, eliminate any that do not and choose one at random from the remaining placements.
\end{description}

In addition we have a rule regarding the initial placement of vertices:

\begin{description}
\item[CP0] If it is possible to place the initial vaccination in time step 0 such that each vaccination is adjacent to at least one other vaccination, then 
eliminate any placements that do not and continue to rule CP1.
\end{description}

\subsection{Hexagonal Grid}
The hexagonal grid is an infinite graph based on a tiling of the plane by regular hexagons. Vertices are vertices of the hexagons and edges are sides. The
hexagon grid is a $3$-regular infinite graph. It is conjectured~\cite{GKP,M} that placing one virus per time step is not enough to contain the outbreak.
On the other hand, it's easy to see that, with two vaccines per time step, the outbreak will be contained after two time steps. 

Let's see how this plays out for the Containment Protocol. Suppose the outbreak starts at vertex $x_0$ at time $t = 0$. Let $N(x_0) = \{x_1,x_2,x_3\}$. 
As we must place vaccines adjacent to vertices in AIV $= \{x_0\}$, the two vaccines at time $t=1$ are two of the three in $N(x_0)$. By symmetry,
all three choices are equivalent. Since there's only one possible placement (up to symmetry) we do not need to consider the five rules, CP0 to CP5, to make a choice.
Let's say that the vaccines are placed at $x_1$ and $x_2$. In time step 1, the fire spreads to $x_3$. Let $N(x_3) = \{x_0, x_4, x_5\}$. Again, there is only
one possible placement of the vaccines, namely, at vertices $x_4$ and $x_5$, which contains the outbreak. Thus, the Containment Protocol matches
the best known containment strategy in the case of the hexagonal grid.

\subsection{Square Grid}

The square grid is an infinite $4$-regular graph based on a tessellation of the plane by squares. It is sometimes called the rectangular or Cartesian grid 
and is the Cartesian product $P_\infty \square P_\infty$ of two infinite paths. It is known~\cite{F,WM} that one vaccine per time step is not enough to 
contain the outbreak. With $n=2$ vaccines per time step, the best strategy requires eight time steps for containment with 18 vertices infected \cite{H,WM}.
We will see that the Containment Protocol matches this optimal performance. 

\begin{figure}[htb]
    \labellist                             
    \small\hair 2pt
    \pinlabel AIV at 63 37
    \pinlabel Vaccination at 78 9
    \pinlabel IV1 at 206 37
    \pinlabel IV2 at 206 9 
    \endlabellist         
\centering
\includegraphics[scale = 1]{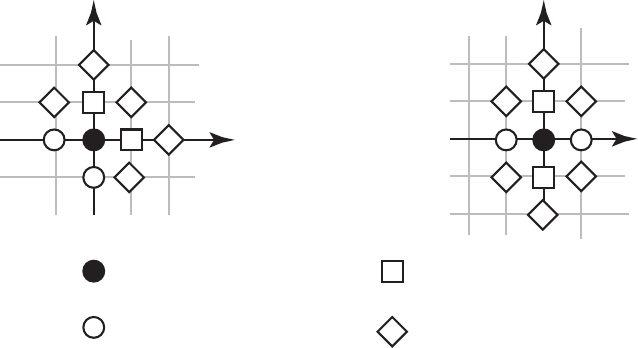}
\caption{Two candidate placements at time 0.}
\label{fig:S0}
\end{figure}
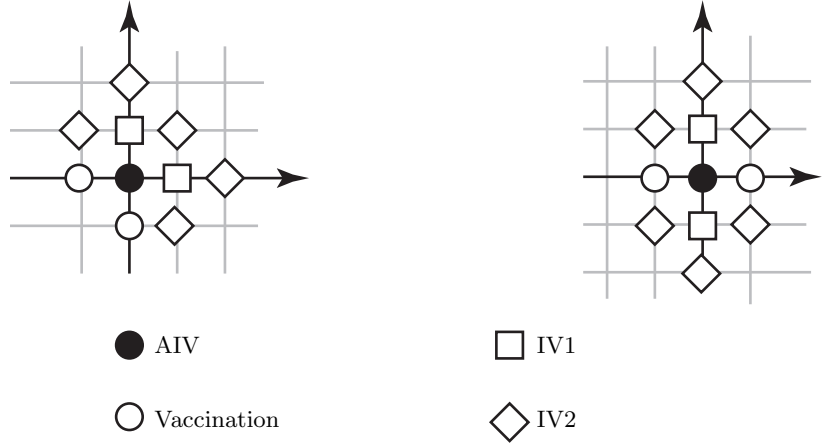

{\bf Time step 0.}
Let's say the outbreak begins at vertex $(0,0)$. Up to symmetry, there are two ways to place two vaccines: 1) at vertices $(-1,0)$ and $(0,-1)$ or 2) at vertices $(-1,0)$ and $(1,0)$. Since neither has the two vaccinations adjacent, following CP0, we continue with both.
As in Figure~\ref{fig:S0}, CP1 says we should choose the first placement for which $|\mbox{IV2}| = 5$, in contrast to $|\mbox{IV2}| = 6$ for the second placement.

 {\bf Time step 1.} Up to symmetry, there are two ways to place vaccines 1) $(-1,1)$ and $(0,2)$; or 2) $(-1,1)$ and $(1,-1)$.
We choose placement 1) which yields $|\mbox{IV2}| = 5$ as opposed to 2) for which $|\mbox{IV2}| = 6$.

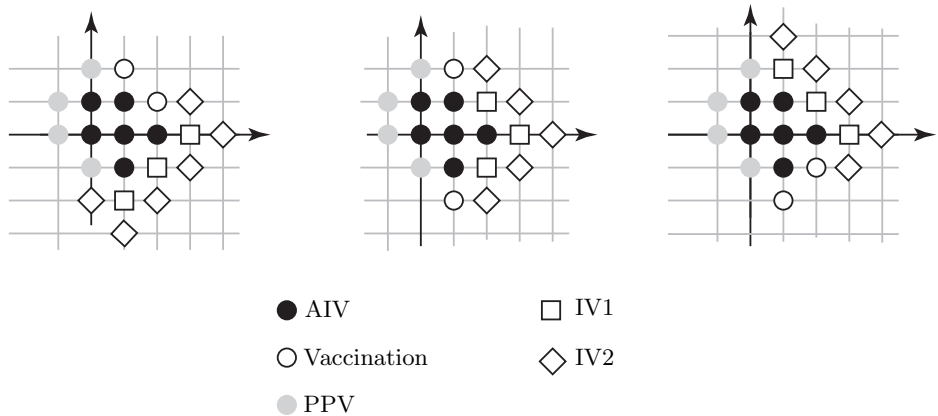
\begin{figure}[htb]
    \labellist                             
    \small\hair 2pt
    \pinlabel AIV at 174 58
    \pinlabel Vaccination at 195 32
    \pinlabel PPV at 175 6
    \pinlabel IV1 at 320 59
    \pinlabel IV2 at 320 32
    \endlabellist         
\centering
\includegraphics[scale = 0.7]{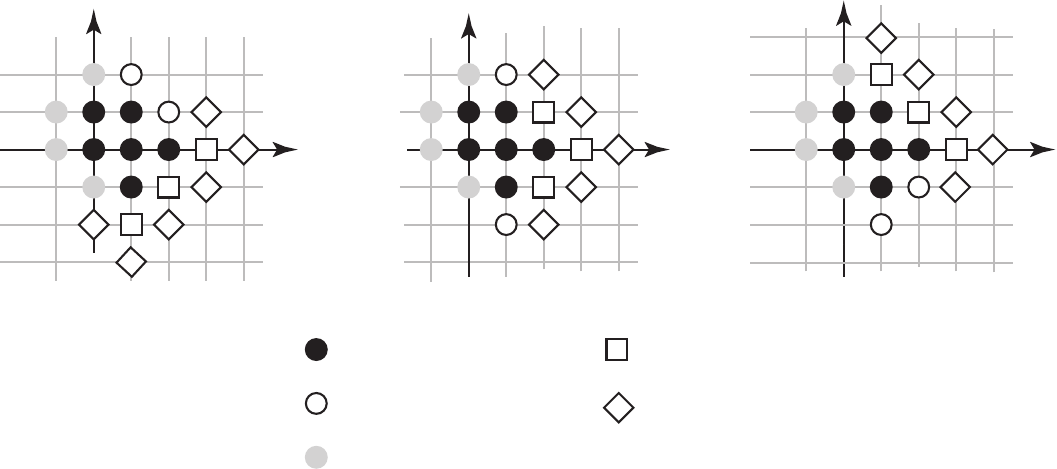}
\caption{Three candidate placements at time 2.}
\label{fig:S2}
\end{figure}

{\bf Time step 2.} There are three possible placements; 1) $(1,2)$ and $(2,1)$ with $|\mbox{IV2}| = 6$; 2)  $(1,2)$ and $(1,-2)$ with $|\mbox{IV2}| = 5$; and
3) $(1,-2)$ and $(2,-1)$ with $|\mbox{IV2}| = 5$, see Figure~\ref{fig:S2}. Using CP1, we eliminate the first placement as a candidate. Neither of the remaining candidates has
adjacent PPV and IV2 vertices, so CP2 does not help us distinguish between them. However, placement 3) has a good PPV at $(0,2)$ (adjacent to the IV1 vertex $(1,2)$ ) while placement 
2) has none. So, by CP3, we use placement 3.

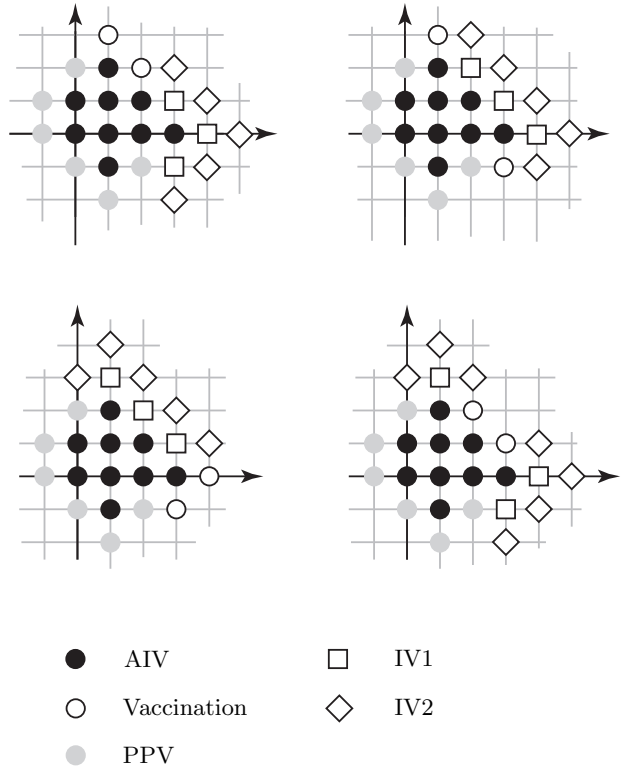
\begin{figure}[htb]
    \labellist                             
    \small\hair 2pt
    \pinlabel AIV at 84 58
    \pinlabel Vaccination at 105 32
    \pinlabel PPV at 85 6
    \pinlabel IV1 at 230 59
    \pinlabel IV2 at 230 32
    \endlabellist       
\centering
\includegraphics[scale = 0.7]{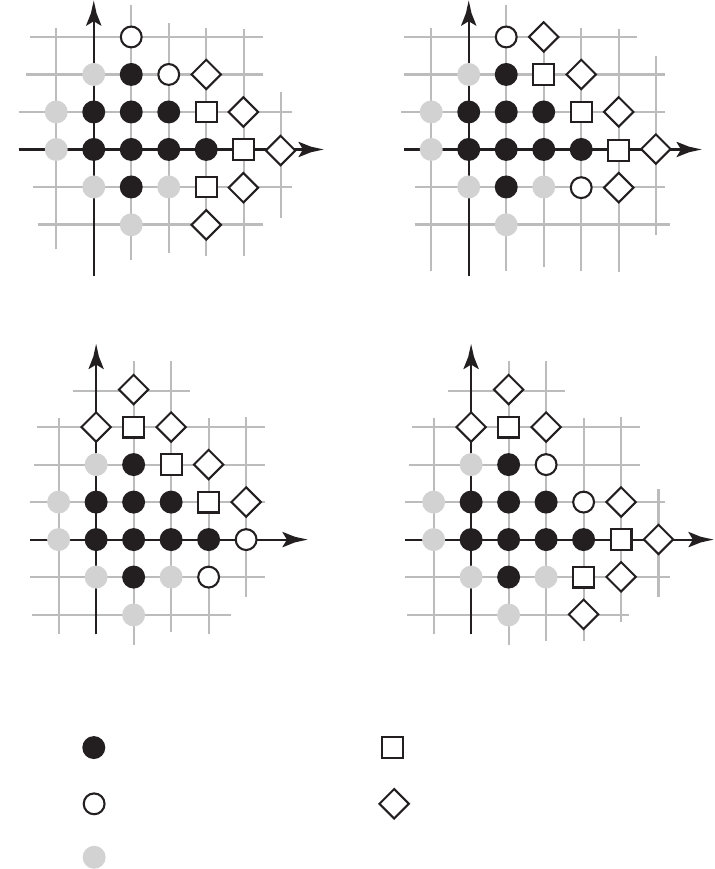}
\caption{Four candidate placements at time 3.}
\label{fig:S3}
\end{figure}

{\bf Time step 3.} The possible placements are 1) $(1,3)$ and $(2,2)$; 2) $(1,3)$ and $(3,-1)$; 3)
$(3,-1)$ and $(4,0)$; and 4) $(2,2)$ and $(3,1)$, see Figure~\ref{fig:S3}. All but the last (with 
$|\mbox{IV2}| = 7$) yield $|\mbox{IV2}| = 5$. By CP1, we eliminate the last placement and 
continue with the first three. Since arrangement 3) has the PPV at $(0,2)$ adjacent to
the IV2 vertex $(0,3)$, we eliminate that one using CP2 and continue to CP3 with the other two. Neither of the remaining placements has a good PPV, 
so we move on to CP4. Since placement 1) has a bad PPV at $(2,-1)$, we will use placement 2). 

{\bf Time step 4.} 1) $(2,3)$ and $(3,2)$ with $| \mbox{IV2}| = 5$; 2) $(2,3)$ and $(4,-1)$ with $| \mbox{IV2}| = 5$; and 3) $(4,-1)$ and $(5,0)$ with
$| \mbox{IV2} | = 4$. By CP1, we use arrangement 3).

{\bf Time step 5.} 1) $(2,4)$ and $(3,3)$ with $| \mbox{IV2}| = 3$; 2) $(2,4)$ and $(5,1)$ with $| \mbox{IV2}| = 3$;  3) $(5,1)$ and $(4,2)$ with
$| \mbox{IV2} |= 4$; and 4) $(3,3)$ and $(4,2)$ with $| \mbox{IV2} |= 5$. By CP1, we eliminate arrangements 3) and 4). CP2 does not allow us to eliminate either of the remaining two candidates. However
CP3 says we should choose arrangement 1) due to the good PPV at $(5,0)$ (adjacent to the IV1 vertex $(5,1)$).

{\bf Time step 6.} 1) $(4,3)$ and $(5,2)$ with $| \mbox{IV2} | = 3$; 2) $(4,3)$ and $(6,1)$ with $| \mbox{IV2} | = 2$; and 3) $(6,1)$ and $(5,2)$ with
$| \mbox{IV2}| = 2$. By CP1, we eliminate arrangement 1). CP2 and CP3 do not allow us to distinguish, but CP4 indicates arrangement 2)
due to the bad PPV at $(3,3)$ in arrangement 3).

\begin{figure}[htb]
\centering
\includegraphics[scale = 1]{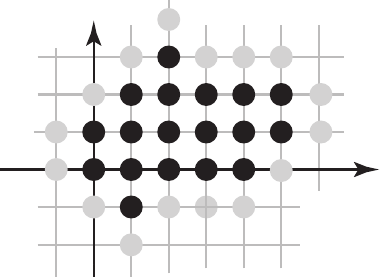}
\caption{The containment protocol stops the outbreak in eight time steps.}
\label{fig:SFinal}
\end{figure}

By Time Step 7, there is only one placement of the two vaccines that satisfy our criteria and that choice ends the outbreak after eight steps and 
with 18 vertices infected, see Figure~\ref{fig:SFinal}. For the square grid, the containment protocol matches the known best solution.

\subsection{Triangular Grid}

\begin{figure}[htb]
\centering
\includegraphics[scale = 0.8]{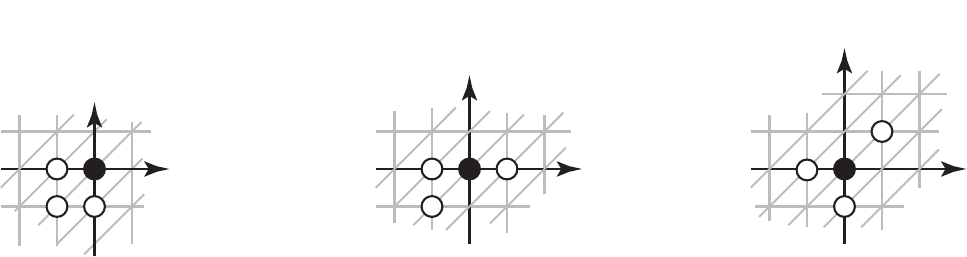}
\caption{Three placements for time step 0.}
\label{fig:T0}
\end{figure}

A triangular grid is the infinite $6$-regular graph that corresponds to a tessellation of the plane by equilateral triangles. For convenience, we place 
this graph in the plane with vertices at lattice points $\{ (m,n) \mid m,n, \in \Z \}$ as in Figure~\ref{fig:T0}. Although Fogarty~\cite{F} proposed an argument showing
that two vaccines per time step do not suffice to contain the outbreak, there is a flaw (see~\cite{GKP}) and this remains a conjecture. It is known that 
the outbreak can be contained by placing three vaccines at each time step. In fact, Messinger~\cite{M} conjectures that the optimal solution sees 17 vertices
infected with containment in six time steps. We now show that the Containment Protocol matches or nearly matches this best known result.

{\bf Time step 0.}
Let's say the outbreak begins at vertex $(0,0)$. Up to symmetry, there are three ways to place three vaccines: 1) at vertices $(-1,0)$, $(-1,-1)$, and $(0,-1)$;
2) at vertices $(-1,0)$, $(-1,-1)$, and $(1,0)$; and 3) at vertices $(-1,0)$, $(1,1)$, and $(0,-1)$,
see Figure~\ref{fig:T0}. By CP0, we use the first placement, the only one that has each of the three vaccinations adjacent to at least one other
vaccination.

{\bf Time step 1.}
Up to symmetry, there are 15 ways to place the three vaccines. Of these, only one, with vaccines at $(-1,1)$, $(0,2)$, and $(1,-1)$, realizes the minimum $| \mbox{IV2} | = 7$.

{\bf Time step 2.}
There are seven vertices adjacent to an AIV, meaning there are $\binom{7}{3} = 35$ placements. However, two of these involve placing a vaccine more than distance
two from the other vacciation sites and PPV vertices, namely 1) $(1,3)$, $(2,-1)$, and $(3,2)$ (for which $(3,2)$ is too far from the other vaccinations) and 2) $(2,-1)$, $(3,0)$, 
and $(3,3)$ (with the last one too distant). Of the remaining placements, only one, at $(1,3)$, $(2,-1)$, and $(3,0)$, realizes the minimum $| \mbox{IV2} | = 6$.

\begin{figure}[htb]
\centering
\includegraphics[scale = 0.9]{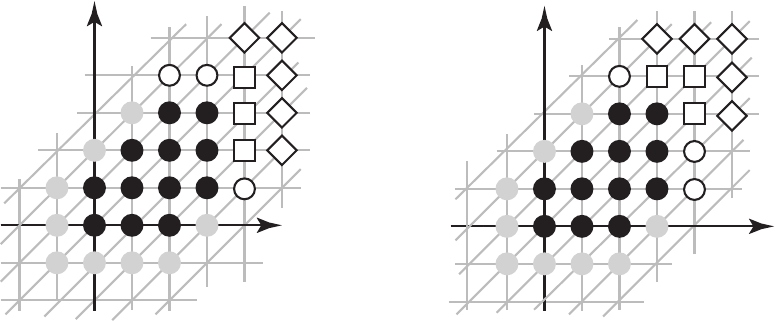}
\caption{Two placements with $| \mbox{IV2}| = 5$ and total distance four.}
\label{fig:T3}
\end{figure}

{\bf Time step 3.} There are six vertices adjacent to an AIV and $\binom{6}{3} = 20$ placements. Four of these realize the minimum $| \mbox{IV2}| = 5$. 
CP2 says we must remove one of the four, with vaccines at $(2,4)$, $(3,4)$, and $(4,4)$, due to the IV2 vertex $(4,0)$ adjacent to the PPV at  $(3,0)$. 
The remaining three all survive CP3 and CP4, and two have the same minimal total distance of four. The one eliminated due to a larger total distance of five, 
has vaccinations at $(2,4)$, $(4,4)$, and $(4,1)$. 

Although the Containment Protocol says we should choose at random between the remaining two placements, we will continue with both to show that 
they both perform well. Going forward we will compare Solution 1, based on vaccinations at $(2,4)$, $(3,4)$, and $(4,1)$, with Solution 2, which uses
the placement $(2,4)$, $(4,1)$, and $(4,2)$.

{\bf Solution 1: Time step 4.}
There are five vertices adjacent to an AIV with $\binom{5}{3} = 10$ placements. Of these, only one realizes the minimum $| \mbox{IV2}| = 3$: 
$(4,5)$, $(5,5)$, and $(5,2)$. 

\begin{figure}[htb]
\centering
\includegraphics[scale = 1]{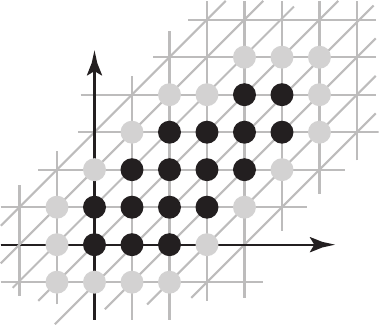}
\caption{Solution one contains the outbreak in six time steps.}
\label{fig:TS1}
\end{figure}

This also means we can contain the outbreak by placing three vaccines at the the three IV2 vertices: $(6,3)$, $(6,4)$, and $(6,5)$.
The complete solution is illustrated in Figure~\ref{fig:TS1}. With this solution, the Containment Protocol matches the best
known result with 17 vertices infected and containment after six time steps.

{\bf Solution 2: Time step 4.}
 As with solution 1, there are five available vertices and ten ways to place vaccines. Five of them realize the minimum $| \mbox{IV2} | = 4$. Two
 are eliminated by CP2. The placement $(3,5)$, $(4,5)$, and $(5,5)$ has an IV2 vertex $(5,2)$ adjacent to the PPV at $(4,2)$. The placement
 $(5,3)$, $(5,4)$, and $(5,5)$ has IV2 $(2,5)$ adjacent to PPV $(2,4)$. The remaining three survive CP3 and CP4 and two realize the 
 minimum total distance four. The one with a larger total distance of five is eliminated by CP5: $(3,5)$, $(5,3)$, and $(5,5)$. Again,
 although the protocol indicates that we should choose one of the remaining two placements at random, we will continue with 
 both to see that they preform equally well. For Solution 2a, we continue with the placement $(3,5)$, $(4,5)$, and $(5,3)$ while
 we will denote the placement $(3,5)$, $(5,3)$, and $(5,4)$ as Solution 2b.
 
\begin{figure}[htb]
\centering
\includegraphics[scale = 0.7]{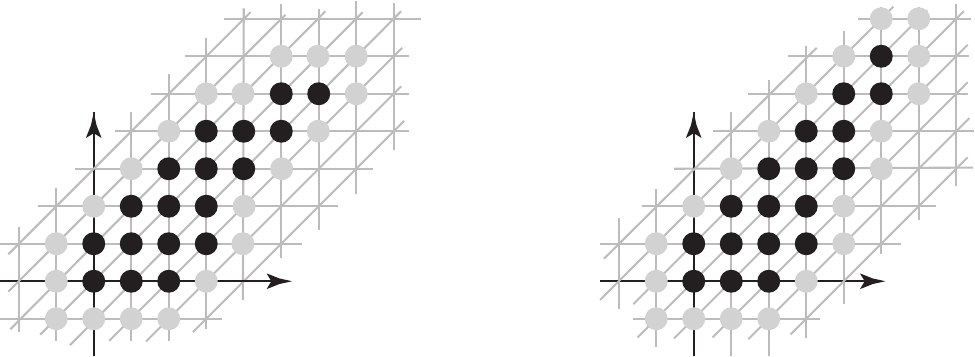}
\caption{Solutions 2a (left) and 2b (right) contain the outbreak in seven time steps.}
\label{fig:TS2}
\end{figure}

 {\bf Solution 2a: Time step 5.}
There are four vertices available, so $\binom{4}{3} = 4$ possible placements. One placement realizes the minimum $| \mbox{IV2} | = 2$:
$(5,6)$, $(6,4)$, and $(6,6)$. This means we can contain the outbreak in one more time step. In fact, we only need two vaccinations
for the last step. In total, 18 vertices are infected, see Figure~\ref{fig:TS2}.

 {\bf Solution 2b: Time step 5.}
Again, there are four vertices and four possible placements. One, $(4,6)$, $(6,5)$, and $(6,6)$, realizes the minimum $| \mbox{IV2} | = 2$.
As in Solution 2a, we contain the outbreak in seven time steps (with only two vaccinations needed in the last step) and 18 vertices
infected, see Figure~\ref{fig:TS2}.

In summary, the Containment Protocol would have us choose at random between three possible solutions. One of these matches the best
known solution with 17 vertices infected and containment in six time steps. The other two are nearly as good with 18 vertices infected
and containment in seven time steps.

\subsection{Strong grid} 
\begin{figure}[htb]
\centering
\includegraphics[scale = 0.7]{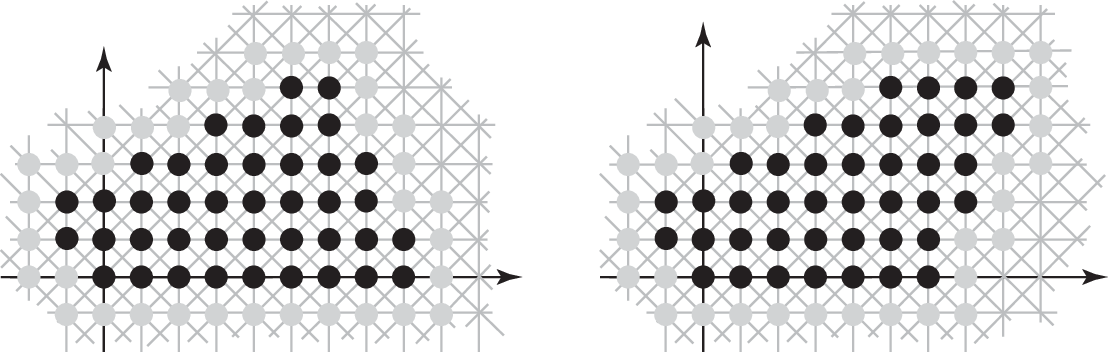}
\caption{Solutions 1a (left) and 1b (right) contain the outbreak in nine time steps.}
\label{fig:SS1}
\end{figure}

The strong grid is an eight regular infinite graph that, similar to the square and triangular grid, we can realize by placing vertices 
at the lattice points $\{ (m,n) \mid m,n, \in \Z \}$. As in Figure~\ref{fig:SS1} each $(m,n)$ is adjacent to its eight nearest neighbors in the plane.
Messinger~\cite{M} shows that four vertices are required to contain an outbreak in the strong grid and conjectures that the 
best solution requires eight time steps. She demonstrates a solution in eight time steps for which 35 vertices are infected.

Let's see how the Containment Protocol handles an outbreak on the strong grid.

{\bf Time step 0.}
Let's say the outbreak begins at vertex $(0,0)$. Up to symmetry, there are five ways to place four vaccines such that each is adjacent to
at least one of the others. Of these, only one results in the minimum $| \mbox{IV2}| = 10$; $(-1,-1)$, $(-1,0)$, $(0,-1)$, and $(1-1)$.
By CP1, we reject the other three : $(-1,-1)$, $(0,-1)$, $(1,0)$, and $(1,1)$ with $|\mbox{IV2}| = 12$;
$(-1,-1)$, $(-1,0)$, $(0,-1)$, and $(1,0)$ with $| \mbox{IV2}| = 13$;
$(-1,-1)$, $(0,-1)$, $(0,1)$, and $(1,0)$ with $| \mbox{IV2}| = 14$;
and $(-1,0)$, $(0,-1)$, $(0,1)$, and $(1,0)$ with $| \mbox{IV2}| = 16$.

{\bf Time step 1.}
There are $\binom{10}{4} = 210$ placements of the four vaccines, but only one realizes the minimum $| \mbox{IV2}| = 10$: $(-2,0)$, $(-2,1)$, 
$(-2,2)$, and $(2,-1)$.

{\bf Time step 2.}
Again, there are $\binom{10}{4} = 210$ placements of the four vaccines. Of these, five realize the minimum $| \mbox{IV2}| = 10$.
Note that two others would also give a low IV2 count, but they are not allowed placements as there is a vaccination more than distance two from the 
others: $(-2,3)$, $(-1,3)$, $(3,-1)$, and $(3,3)$; and $(-2,3)$, $(3,-1)$, $(3,0)$, and $(3,3)$, both of which have
$(3,3)$ too distant from the other vaccinations and the PPV set. 
The five placements are not distinguished by CP2, but four have a good PPV. The one that does not is $(-2,3)$, $(-1,3)$, $(3,-1)$,
and $(3,0)$. The remaining four are not distinguished by CP4. Two of them minimize the total distance as five. The two that do not are
$(-2,3)$, $(2,3)$, $(3,-1)$, and $(3,3)$; and $(-2,3)$, $(3,-1)$, $(3,2)$, and $(3,3)$. As usual, although the Containment Protocol requires us to 
choose one of the remaining two placements at random, we will instead continue with both to see how the protocol performs. For
Solution 1: $(-2,3)$, $(-1,3)$, $(0,3)$, and $(3,-1)$; and Solution 2: $(-2,3)$, $(3,-1)$, $(3,0)$, and $(3,1)$. 

{\bf Solution 1: Time step 3.} 
Again there are $\binom{10}{4} = 210$ placements of the four vaccines and five realize the minimum
$|\mbox{IV2}| = 10$. Two others would give a low count but are not allowed as there is a vaccination too far from the others: $(0,4)$, $(1,4)$,
$(4,-1)$, and $(4,4)$ where $(4,4)$ is too far; and $(0,4)$, $(4,-1)$, $(4,0)$, and $(4,4)$, again with $(4,4)$ too far. The remaining five placements
are not distinguished by CP2. Three have a good PPV and the other two, which are eliminated, are $(0,4)$, $(1,4)$, $(4,-1)$, and $(4,0)$; and
$(0,4)$, $(4,-1)$, $(4,0)$,  and $(4,1)$. Of the remaining three, two have a bad PPV and are eliminated by CP4: $(0,4)$, $(3,4)$, $(4,-1)$, and
$(4,4)$; and $(0,4)$, $(4,-1)$, $(4,3)$, and $(4,4)$, both having a good PPV at $(3,-1)$ and a bad one at $(0,3)$.  The remaining placement is the one we will continue with; $(0,4)$, $(1,4)$, $(2,4)$, and $(4,-1)$ with a good PPV at $(3,-1)$.

{\bf Solution 1: Time step 4.} 
Once again, there are $\binom{10}{4} = 210$ placements of the four vaccines and seven realize the minimum $| \mbox{IV2}| = 10$. 
Another, $(2,5)$, $(5,-1)$, $(5,0)$ and $(5,5)$ would also give a low count but the vaccination at $(5,5)$ is more than distance two from the others.
We reject one of the seven using CP2: $(2,5)$, $(3,5)$, $(4,5)$, and $(5,5)$ due to the IV2 vertex at $(4,-2)$ which is adjacent to the PPV $(4,-1)$.
CP3 allows us to eliminate two placements that have no good PPV: $(2,5)$, $(3,5)$, $(5,-1)$, and $(5,0)$; and $(2,5)$, $(5,-1)$, $(5,0)$, and $(5,1)$.
Two more are removed due to a bad PPV: $(2,5)$, $(4,5)$, $(5,5)$, and $(5,-1)$; and $(2,5)$, $(5,-1)$, $(5,4)$, and $(5,5)$, both having a good
PPV at $(4,-1)$ and a bad one at $(2,4)$. The remaining two placements are distinguished by CP5. The minimum total distance is five for $(2,5)$,
$(3,5)$, $(4,5)$, and $(5,-1)$ and we will continue with this placement. We reject the placement with total distance six: $(2,5)$, $(3,5)$, $(5,-1)$, and
$(5,5)$.

{\bf Solution 1: Time step 5.} We continue to have $\binom{10}{4} = 210$ placements. In this case, there is a unique placement that realizes the
minimum $| \mbox{IV2}| = 8$: $(4,6)$, $(5,6)$, $(6,6)$, and $(6,-1)$.

{\bf Solution 1: Time step 6.} Now there are $\binom{8}{4} = 70$ placements.  Of these, three realize the minimum $| \mbox{IV2}| = 6$ and
none of those have a PPV adjacent to an IV2.  Two have a good PPV. By CP3, we eliminate the third: $(7,-1)$, $(7,0)$, $(7,5)$, and $(7,6)$.
The remaining two have no bad PPV and both have a total distance of five. Ordinarily, the Containment Protocol would have us choose one of them
at random. Instead we will continue with both. For Solution 1a: $(7,-1)$, $(7,4)$, $(7,5)$, and $(7,6)$ with a good PPV at $(6,-1)$; 
for Solution 1b: $(7,-1)$, $(7,0)$, $(7,1)$, and $(7,6)$ with a good PPV at $(6,6)$.

{\bf Solution 1a: Time step 7.}
We are down to $\binom{6}{4} = 15$ placements. There are three that realize the minimum $| \mbox{IV2}| = 4$. This means all three are complete
in one more time step with the same amount of vertices infected. There is one that we eliminate by CP3, as it's the only one lacking a good PPV:
$(8,-1)$, $(8,0)$, $(8,3)$, and $(8,4)$. From the remaining two, we randomly chose the first: $(8,-1)$, $(8,2)$, $(8,3)$, $(8,4)$ with a good PPV
at $(7,-1)$, rather than the second: $(8,-1)$, $(8,0)$, $(8,1)$, and $(8,4)$ with a good PPV at $(7,4)$. 

This leads to the solution shown at left in Figure~\ref{fig:SS1}; the outbreak is contained in nine time steps with 41 vertices infected.

{\bf Solution 1b: Time step 7.} 
We are down to $\binom{6}{4} = 15$ placements. Of these, three realize the minimum $| \mbox{IV2}| = 4$. One lacks a good PPV
and is rejected by CP3: $(8,1)$, $(8,2)$, $(8,5)$, and $(8,6)$. The other two are equal with respect to CP4 and CP5 and also in terms of the number of vertices infected 
and time steps for the final solution. So, we choose one at random: $(8,1)$, $(8,2)$, $(8,3)$, and $(8,6)$ with a good PPV at $(7,6)$; and reject the other
$(8,1)$, $(8,4)$, $(8,5)$, and $(8,6)$ with a good PPV at $(7,1)$.

This leads to the solution shown at right in Figure~\ref{fig:SS1}; again, the outbreak is contained in nine time steps with 41 vertices infected.

{\bf Solution 2: Time step 3.} 
There are $\binom{10}{4} = 210$ placements of which eight would yield $| \mbox{IV2}| = 10$. One of these, $(-2,4)$, $(-1,4)$, $(4,1)$, and $(4,4)$, has 
$(4,4)$ more than two away from the other vaccination sites and is not allowed for that reason. Another is rejected due to CP2: $(4,1)$, $(4,2)$, $(4,3)$, and $(4,4)$,
since the IV2 $(-3,3)$ is adjacent to the PPV $(-2,3)$. Next, there are placements that have a good PPV, so we reject the two that do not: $(-2,4)$, $(-1,4)$,
$(0,4)$, and $(4,1)$; and $(-2,4)$, $(-1,4)$, $(4,1)$, and $(4,2)$. Of the remaining four placements, we reject two that have a  bad PPV: $(-2,4)$, $(3,4)$,
$(4,1)$, and $(4,4)$; and $(-2,4)$, $(4,1)$, $(4,3)$, and $(4,4)$, both of which have $(3,1)$ as a bad PPV.
 
The remaining two placements are distinguished by total distance (see CP5). The minimal distance is five, which is achieved by $(-2,4)$, $(4,1)$,
$(4,2)$, and $(4,3)$, with which we continue. The final placement $(-2,4)$, $(4,1)$, $(4,2)$, and $(4,4)$ is rejected due to a total distance of six.

{\bf Solution 2: Time step 4.} There are $\binom{10}{4} = 210$ placements only one of which realizes the minimum $| \mbox{IV2}| = 8$: $(-2,5)$,
$(5,3)$, $(5,4)$, and $(5,5)$.

{\bf Solution 2: Time step 5.} 
There are $\binom{8}{4} = 70$ placements and three realize the minimum $| \mbox{IV2}| = 6$. Two have a good PPV, so, by CP3, we reject the 
third: $(-2,6)$, $(-1,6)$, $(4,6)$, and $(5,6)$. The remaining two placements are not distinguished by CP4 or CP5. Ordinarily, the Containment
Protocol would have us choose between them at random, but we will carry both forward to see how the protocol performs in all cases.
For solution 2a, we place vaccines at $(-2,6)$, $(-1,6)$, $(0,6)$, and $(5,6)$. In this case, $(5,5)$ is a good PPV.
For solution 2b, we continue with $(-2,6)$, $(3,6)$, $(4,6)$, and $(5,6)$ which has $(-2,5)$ as a good PPV. 

\begin{figure}[htb]
\centering
\includegraphics[scale = 0.9]{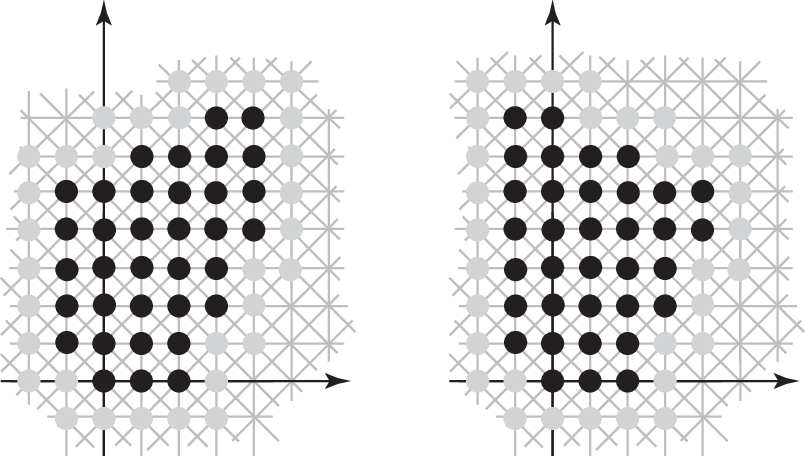}
\caption{Solutions 2a (left) and 2b (right) contain the outbreak in eight time steps.}
\label{fig:SS2}
\end{figure}

{\bf Solution 2a: Time step 6.}
There are $\binom{6}{4} =  15$ placements and two realize the minimum $| \mbox{IV2}| = 4$. One has a good PPV, so, by CP3, we reject the other:
$(0,7)$, $(3,7)$, $(4,7)$, and $(5,7)$ and continue with $(0,7)$, $(1,7)$, $(2,7)$, and $(5,7)$. In the next time step, we place vaccines at the four
IV2 vertices and contain the outbreak in eight time steps with 35 vertices infected, see Figure~\ref{fig:SS2}.

{\bf Solution 2b: Time step 6.}
There are $\binom{6}{4} =  15$ placements and three realize the minimum $| \mbox{IV2}| = 4$. Of these, only one has a good PPV, so by CP3, we 
reject the other two: $(-2,7)$, $(-1,7)$, $(0,7)$, and $(3,7)$; and $(-2,7)$, $(-1,7)$, $(2,7)$, and $(3,7)$. We continue with the placement
$(-2,7)$, $(1,7)$, $(2,7)$, and $(3,7)$. We next place vaccines at the four IV2 vertices to terminate the outbreak in eight time steps, with
35 vertices infected, see Figure~\ref{fig:SS2}.

In summary, the Containment Protocol produces one of four solutions, all with equal likelihood. Two of the four match the best known solution with 
containment in eight time steps and 35 vertices infected. The other two are close to optimal with 41 vertices infected after nine time steps.

\section{
\label{sec:PG}%
The pentagon graph}

\begin{figure}[htb]
\centering
\includegraphics[scale = 1]{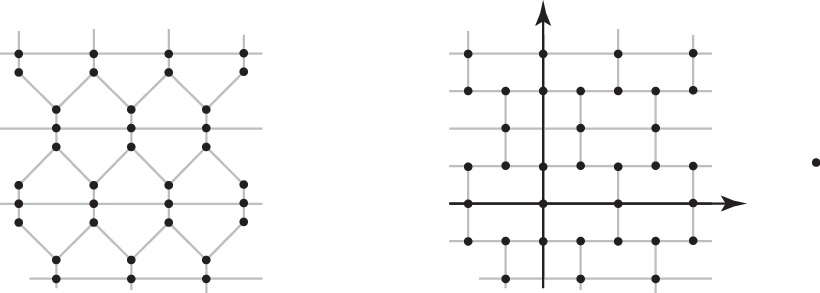}
\caption{A pentagonal tiling (left). The resulting graph on the lattice (right).}
\label{fig:pent}
\end{figure}

We have seen that the Containment Protocol performs quite well, often meeting the best known outcome for some well-studied regular infinite graphs.
One virtue of the protocol is that it can be applied to virus spread on any graph. As an example, we
investigate a graph resulting from a pentagonal tiling of the plane as in Figure~\ref{fig:pent}. This results in an infinite
graph that is not regular as it has vertices of degree three and four. 

\begin{conjecture}
\label{conj:pent}
One vaccination per time step is not sufficient to control an outbreak on the Pentagon Graph.
\end{conjecture}

It is easy to see that, much like the hexagonal grid, an outbreak that begins on a degree three vertex can be contained in two time 
steps if there are two vaccinations made per time step. Let's see how the Containment Protocol handles this situation. As in Figure~\ref{fig:pent},
we imagine the outbreak beginning at $(0,1)$. Up to symmetry there are two ways to place the two vaccinations: $(-1,1)$ and $(0,0)$ with
$| \mbox{IV2}| = 2$; and $(-1,1)$ and $(1,1)$ with $| \mbox{IV2}| = 3$. Both placements fail to have the vaccinations adjacent (see CP0), and
CP1 says we should continue with the first placement. We next place our two vaccinations at the two IV2 vertices to contain the outbreak
in two steps with two vertices infected, which is the optimal solution.

We next consider an outbreak that starts at the degree four vertex $(0,0)$ and allow ourselves two vaccinations per time step.

{\bf Time step 0.} Up to symmetry, there are three ways to place two vaccinations adjacent to the initial degree four vertex at $(0,0)$. For none
of those three are the two vaccinations adjacent, so CP0 does not apply.  However, only one 
of them realizes the minimum $|\mbox{IV2}| = 4$: $(-2,0)$ and $(2,0)$.

{\bf Time step 1.} Up to symmetry, there are three ways to place the next two vaccinations, none of which can be distinguished by CP1 to CP5.
Ordinarily, the Containment Protocol would have us choose one at random, but we will continue our practice of following through with 
all three: Solution 1: $(-1,1)$ and $(1,1)$; Solution 2: $(-1,-1)$ and $(-1,1)$; and Solution 3: $(-1,1)$ and $(1,-1)$. 

{\bf Solution 1: Time step 2.} Among the $\binom{4}{2} = 6$ ways to place two vaccinations, the only one that realizes the minimum 
$|\mbox{IV2}| = 2$ is $(-1,-2)$ and $(1,-2)$. Placing two vaccinations at the IV2 positions contains the outbreak in four time steps with
seven vertices infected.

\begin{figure}[htb]
    \labellist                             
    \small\hair 2pt
    \pinlabel $0$ at 28 43
    \pinlabel $0$ at 100 43
    \pinlabel $1$ at 45 61
    \pinlabel $1$ at 82 61
    \pinlabel $2$ at 45 8
    \pinlabel $2$ at 82 8
    \pinlabel $3$ at 10 25
    \pinlabel $3$ at 118 25
    \pinlabel $0$ at 208 43
    \pinlabel $0$ at 280 43
    \pinlabel $1$ at 226 61
    \pinlabel $1$ at 227 25
    \pinlabel $2$ at 261 80
    \pinlabel $2$ at 261 8
    \pinlabel $3$ at 299 61
    \pinlabel $3$ at 299 25
    \pinlabel $0$ at 388 43
    \pinlabel $0$ at 460 43
    \pinlabel $1$ at 406 61
    \pinlabel $1$ at 442 25
    \pinlabel $2$ at 406 7
    \pinlabel $2$ at 442 79
    \pinlabel $3$ at 369 25
    \pinlabel $3$ at 478 61
    \endlabellist       
\centering
\includegraphics[scale = 0.7]{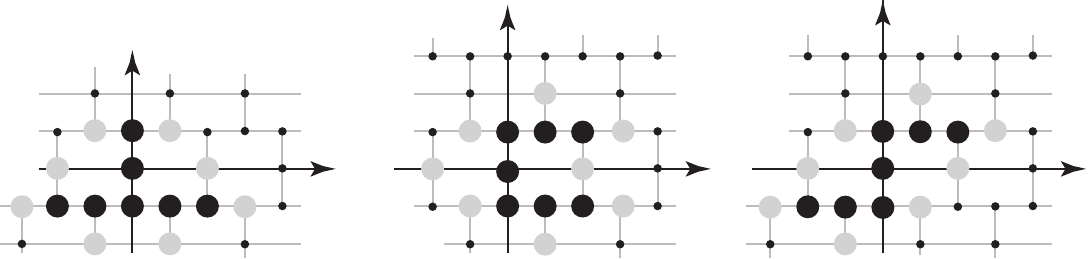}
\caption{Starting from a single degree four vertex, Solution 1 (left), Solution 2 (center), and Solution 3 (right).
\label{fig:P1}
}
\end{figure}
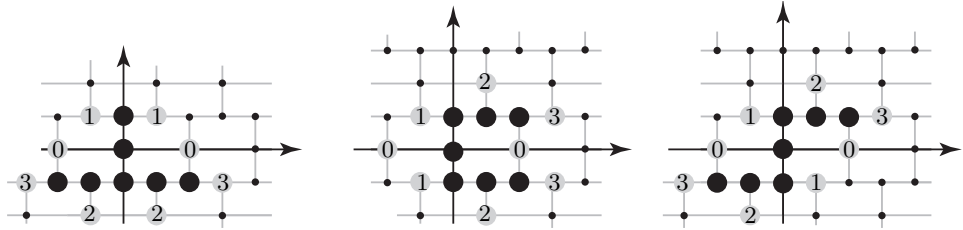

Solutions 2 and 3 are similar, also stopping the outbreak in four time steps with seven vertices infected, see Figure~\ref{fig:P1}.

As the analysis of an outbreak at a single vertex was straightforward, we also investigated what would happen if an outbreak began with a 
pair of adjacent vertices infected. There are three cases as the two vertices may have the same degree three or four, or different degrees.

{\bf Two degree three vertices infected:}
Suppose an outbreak begins with two adjacent degree three vertices infected at $(0,1)$ and $(1,1)$.
Given Conjecture~\ref{conj:pent}, it is likely not possible
to contain the outbreak with a single vaccine per time step, so we will allow ourselves two. 

{\bf Time step 0.}
Up to symmetry, there are four ways to place two vaccinations adjacent to the infected vertices. None have the vaccinations adjacent to one another, 
so CP0 does not help. Two of the four realize the minimum $| \mbox{IV2}| = 4$, and we will continue with both, even though 
the Containment Protocol says we should pick one at random. Solution 1: $(0,0)$ and $(1,2)$ or Solution 2: $(-1,1)$ and $(1,2)$.

{\bf Solution 1: Time step 1.} Up to symmetry, there are four ways to place the two vaccinations and all four result in $|\mbox{IV2}| = 4$. 
Using CP2, we eliminate all but one placement: $(-2,1)$ and $(3,1)$.

{\bf Solution 1: Time step 2.} Up to symmetry, there are four ways to place the two vaccinations and only one realizes the minimum 
$|\mbox{IV2}| = 4$: $(-3,2)$ and $(4,0)$. 

{\bf Solution 1: Time step 3.} Up to symmetry, there are four ways to place the two vaccinations. These four are all equivalent with respect to 
CP1 through 5. The Containment Protocol would have us choose one at random. Instead we will continue with all four so that we can compare
the results: Solution 1a: $(-2,3)$ and $(0,3)$; Solution 1b: $(-2,3)$ and $(1,-1)$; Solution 1c: $(-2,3)$ and $(3,-1)$; and Solution 1d:
$(0,3)$ and $(1,-1)$.

\begin{figure}[htb]
\centering
\includegraphics[scale = 0.9]{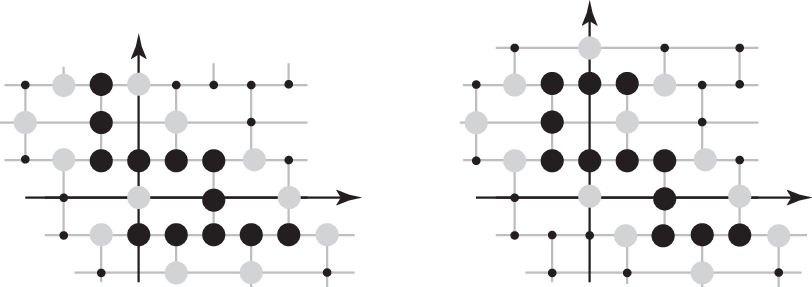}
\caption{Solution 1a (left) and 1b (right).}
\label{fig:PTTS11}
\end{figure}

{\bf Solution 1a: Time step 4.} There are $\binom{4}{2} = 6$ ways to place vaccines and only one realizes the minimum $| \mbox{IV2}| = 2$:
$(1,-2)$ and $(3,-2)$. Covering those two vertices in the next step contains the outbreak in six time steps with 12 vertices infected,
see Figure~\ref{fig:PTTS11}.

{\bf Solution 1b: Time step 4.} Again, there are $\binom{4}{2} = 6$ ways to place vaccines and only one realizes the minimum $| \mbox{IV2}| = 2$:
$(0,4)$ and $(3,-2)$. Placing vaccines at those two vertices contains the outbreak in six time steps with 12 vertices infected,
see Figure~\ref{fig:PTTS11}.

\begin{figure}[htb]
\centering
\includegraphics[scale = 0.9]{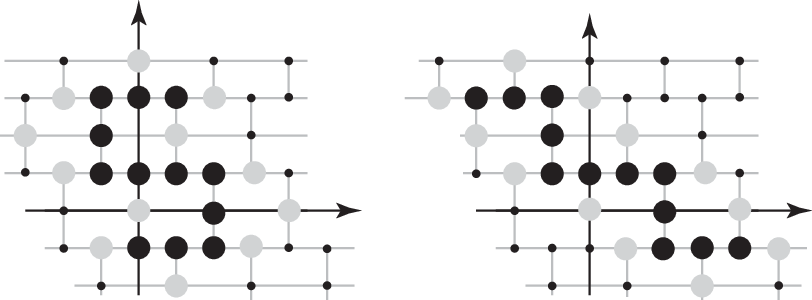}
\caption{Solution 1c (left) and 1d (right).}
\label{fig:PTTS12}
\end{figure}

{\bf Solution 1c: Time step 4.} There are $\binom{4}{2} = 6$ ways to place vaccines and only one realizes the minimum $| \mbox{IV2}| = 2$:
$(0,4)$ and $(1,-2)$. Covering those two vertices in the next step contains the outbreak in six time steps with 12 vertices infected,
see Figure~\ref{fig:PTTS12}.

{\bf Solution 1d: Time step 4.} Again, there are $\binom{4}{2} = 6$ ways to place vaccines and only one realizes the minimum $| \mbox{IV2}| = 2$:
$(-2,4)$ and $(3,-2)$. Covering those two vertices in the next step contains the outbreak in six time steps with 12 vertices infected,
see Figure~\ref{fig:PTTS12}.

{\bf Solution 2: Time step 1.} There are $\binom{4}{2} = 6$ ways to place vaccines and three realize the minimum $| \mbox{IV2}| = 4$.
However, only one has no PPV adjacent to an IV2, so we continue with that one: $(-2,0)$ and $(3,1)$.

\begin{figure}[htb]
\centering
\includegraphics[scale = 1]{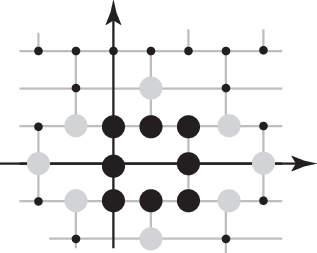}
\caption{Solution 2.}
\label{fig:PTTS2}
\end{figure}

{\bf Solution 2: Time step 2.} There are $\binom{4}{2} = 6$ ways to place vaccines but only one realizes the minimum $| \mbox{IV2}| = 2$:
$(-1,-1)$ and $(4,0)$. Then placing vaccines at the two IV2 vertices completes the outbreak in four time steps with eight vertices infected,
see Figure~\ref{fig:PTTS2}

In summary, for an outbreak that begins with two adjacent degree three vertices infected, using the Containment Protocol, we found 
two types of solutions. Solutions 1a, 1b, 1c, and 1d contain the outbreak in six time steps with 12 vertices infected. The best outcome is
solution 2 with eight vertices infected after four time steps. 

{\bf Two vertices including one of degree four infected:}
If an infection begins with two adjacent vertices infected, one of degree four, the analysis of the Containment Protocol becomes quite involved,
branching off into many possible solutions with a wide range of time steps and number of vertices infected. Rather than chasing through all the details,
we will present the best solutions we found using the Containment Protocol. We conjecture that these are the best possible solutions.

\begin{figure}[htb]
    \labellist                             
    \small\hair 2pt
    \pinlabel $0$ at 28 64
    \pinlabel $0$ at 46 82
    \pinlabel $1$ at 82 101
    \pinlabel $1$ at 136 64
    \pinlabel $2$ at 118 82
    \pinlabel $2$ at 46 28
    \pinlabel $3$ at 10 46
    \pinlabel $3$ at 82 11
    \pinlabel $4$ at 118 11
    \pinlabel $4$ at 154 28
    \pinlabel $5$ at 173 46        
    \pinlabel $0$ at 245 64
    \pinlabel $0$ at 353 64
    \pinlabel $1$ at 263 46
    \pinlabel $1$ at 337 46
    \pinlabel $2$ at 299 28
    \pinlabel $2$ at 263 101
    \pinlabel $3$ at 227 82
    \pinlabel $3$ at 371 101
    \pinlabel $4$ at 389 82
    \pinlabel $4$ at 281 119
    \pinlabel $5$ at 317 137
    \pinlabel $5$ at 353 137
    \pinlabel $6$ at 389 119
    \endlabellist       
\centering
\includegraphics[scale = 0.9]{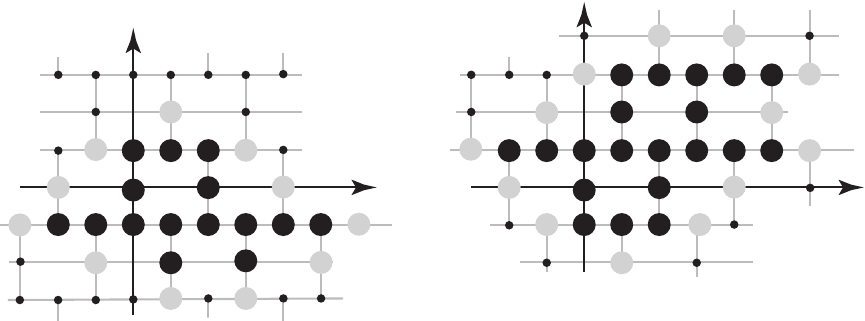}
\caption{Conjectural best solutions for an outbreak that starts with two adjacent vertices, one of degree 4.}
\label{fig:PFFS}
\end{figure}
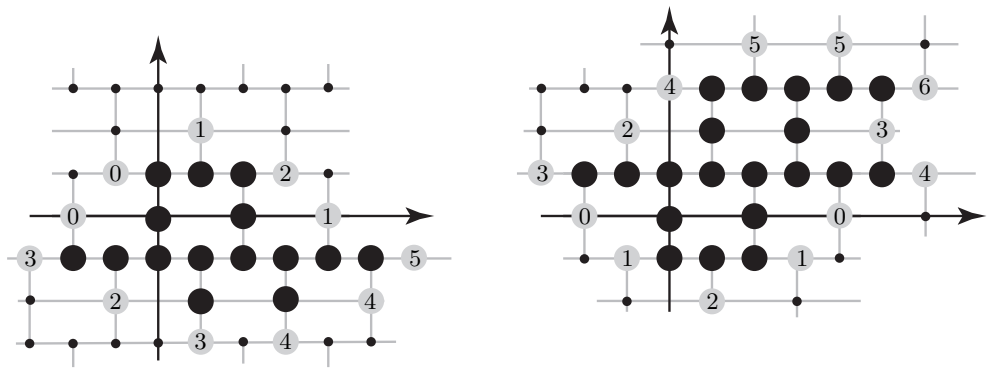

There are two cases. Suppose first that an outbreak begins with the adjacent vertices $(0,0)$, of degree four, and $(0,1)$, of degree three. 
Using the Containment Protocol, the best solution that we found is at left in Figure~\ref{fig:PFFS} and entails containing the outbreak in six time steps with 15 vertices infected. Next, imagine an outbreak that starts with vertices $(0,0)$ and $(2,0)$, both of degree four. 
To the right in Figure~\ref{fig:PFFS}, we present the best solution we found with the Containment Protocol, which sees the outbreak terminated in
seven times steps with 20 vertices infected.

\end{document}